\documentstyle[11pt]{article}

\normalbaselines

\input{mssymb}

\newtheorem{theorem}{Theorem}[section]

\newtheorem{definition}[theorem]{Definition}

\def\A{A^{ij}_k}

\def\C{{\Bbb C}}

\def\E{{\cal E}}

\def\Z{{\Bbb Z}}
\def\g{\hbox{${g}$}}
\def\Ug{\hbox{${\cal U}_g$}}
\def\E{\hbox{${\cal E}$}}

\begin{document}
\bibliographystyle{unsrt}
\vbox{\vspace{38mm}}
\begin{center}
{\LARGE \bf Differential Calculi of Poincar\'e-Birkhoff-Witt type\\[2mm]
on Universal Enveloping Algebras}\\[5mm]

R. Martini, G.F. Post and P.H.M. Kersten\\
{\it University of Twente\\ Department of Applied Mathematics\\
P.O. Box 217, 7500 AE Enschede\\
The Netherlands}\\[5mm]
April 12, 1995
\end{center}

\begin{abstract}
Differential calculi of Poincar\'e-Birkhoff-Witt type on universal
enveloping algebras of Lie algebras \g\ are defined. This definition
turns out to be independent of the basis chosen in \g. The role of
automorphisms of \g\ is explained.  It is proved that no differential
calculus of Poincar\'e-Birkhoff-Witt type exists on semi-simple Lie
algebras. Examples are given, namely $gl_n$, Abelian Lie algebras, the
Heisenberg algebra, the Witt and the Virasoro algebra. Completely
treated are the 2-dimensional solvable Lie algebra, and the
3-dimensional Heisenberg algebra.
\\[2mm]
A.M.S. Subject Classification (1991): 58A10, 16S30, 81R50.\\[1mm]
PACS: 210, 240, 0365.
\end{abstract}

\section{Introduction}

Recently there has been considerably interest in non-commutative
differential geometry, both mathematically and as a framework for
certain physical theories. In particular there is much activity in
differential geometry related to quantum groups.

For instance Woronowicz \cite{Wor} developed a general theory of
fields and forms on quantum groups. This general theory has been
reformulated by Wess and Zumino \cite{WZ} in an easier accessible way.
Furthermore their approach may turn out to be more suitable for
explicit computations with regard to applications in physics. A number
of very interesting papers (e.g. \cite{SWZ}-\cite{M}) elucidating
various aspects and studying specific examples have been written
since.

In this paper we deal with possible differential calculi on Lie
algebras. To be more precise, differential calculi on the universal
enveloping algebra of a given Lie algebra are defined and we
investigate the possibility to construct such a calculus. Our
approach, starting from generators and relations, is somewhat related
to the method of Wess and Zumino, although it is different from the
traditional approaches to construct differential calculi on quantum
groups.  For a Hopf algebraic interpretation of the method and
results, see
\cite{MH}.

It can be proved that there exists no differential calculus in this
sense on any semisimple Lie algebra. However for non-semisimple Lie
algebras things turn out differently. For example there exists a
number of interesting non-trivial differential calculi on any general
linear algebra, the Heisenberg algebra and also on the
infinite-dimensional Virasoro algebra.

\section{Conditions for differential calculi of
Poincar\'e-Birkhoff-Witt type on the universal
enveloping algebra of a Lie algebra}
\setcounter{equation}{0}

In this section our aim is to determine conditions for differential
calculi of Poincar\'e-Birkhoff-Witt type on the universal enveloping
algebra of a given Lie algebra \g. Let us start to describe what we
mean by this. If $x^i,\ (i\in I)$ is a basis of \g, and $\preccurlyeq$
is a well-ordering on $I$, then the Poincar\'e-Birkhoff-Witt theorem
tells that
\[ x^{i_1}x^{i_2}\dots x^{i_m} \qquad  i_1 \preccurlyeq i_2
\preccurlyeq \dots \preccurlyeq i_m \]
is a basis of the universal enveloping algebra \Ug.  We shall use the
following
\begin{definition}
By a differential calculus of Poincar\'e-Birkhoff-Witt type we mean an
algebra extension \E\ of \Ug, such that
\begin{enumerate}
\item[(a)] $ \E = \oplus_{n=0}^{\infty} \E^n$  is a grading of \E, with
$\E^0 = \Ug$.
\item[(b)] There is a linear map d: $\E \to \E$, d$(\E^n) \subset
\E^{n+1}$ that satisfies
\begin{itemize}
  \item[(1)] $\mbox{d}^2 = 0.$
  \item[(2)] d$(\omega \eta)= \mbox{d}(\omega)\eta + (-1)^n\omega
             \mbox{d}\eta$ for $\omega \in \E^n$.
\end{itemize}
\item[(c)] If $y^i=d(x^i)$ then the elements
\[ y^{j_1}y^{j_2}\dots y^{j_n}x^{i_1}x^{i_2}\dots x^{i_m} \]
with $j_1 \prec \dots \prec j_n$ and $i_1 \preccurlyeq i_2
\preccurlyeq \dots \preccurlyeq i_m$ form a (linear) basis of $\E^n$.
\end{enumerate}
\end{definition}
Due to these requirements we have the following relations in \E.
First we have
\begin{equation}
x^i x^j =x^j x^i + c^{ij}_k x^k \qquad (j \prec i)
\label{2.1}
\end{equation}
where $c^{ij}_k$ are the structure constants of the given Lie algebra
\g. Throughout the summation convention is applied.

Secondly we have
\begin{equation}
x^i y^j = y^j x^i + y^k A^{ij}_k \qquad (i,j \in I)
\label{2.2}
\end{equation}

We assume that $A^{ij}_k$ are (complex) numbers. In the case that all
$A^{ij}_k =0$, the relations (\ref{2.2}) represent the commutation
relations of the standard calculus. Our task will be to derive
conditions for the $\A$ (expressed in terms of $c^{ij}_k$) in order to
get a calculus of Poincar\'e-Birkhoff-Witt type. One might object that
the requirement that the $A^{ij}_k$ are constant is rather severe.
However, one can prove that adding non-zero terms, linear in $x^i$
makes it impossible to define a Hopf algebra structure on \E, see
\cite{MH}.  This aspect we will discuss here any further.

Finally we have, by applying d to the relations (\ref{2.2}), and using
the requirements from definition 2.1 (b)
\begin{equation}
y^i y^j = -y^j y^i \quad ( j \prec i) \quad \mbox{ and }\quad y^iy^i=0
\end{equation}

We introduce an ordering on the monomials $z^{i_1}z^{i_2}\dots
z^{i_m}$, where $z$ is either $x$ or $y$. First we order by length,
and after that lexicographically, with $x^i$ larger then $y^j$. Doing
this the relations (\ref{2.1}), (\ref{2.2}) and (2.3) can be seen as
rewrite-rules, replacing the leading term (in the left-hand side) by
lower terms. With respect to these rewrite-rules we have the ``normal
monomials'', i.e. the monomials that can not be simplified using the
rewrite-rules (\ref{2.1}), (\ref{2.2}) and (2.3). It is easy to see
that the elements in definition 2.1 (c) are exactly the normal
monomials.  Now we can apply the diamond lemma (see i.e. \cite{B}). It
states the following: if rewriting $(z^iz^j)z^k$ and $z^i(z^jz^k)$
lead to the same normal form, then the normal monomials are a basis of
\E. Here $z=x$ or $z=y$ with $z^k \prec z^j \prec z^i$. So we have the
following steps:
\begin{enumerate}
\item[{\bf I.}] Consider $(x^ix^j)x^k$ and $x^i(x^jx^k)$. In this case
Jacobi's identity guarantees the same normal forms. (It is in fact a
way to prove the Poincar\'e-Birkhoff-Witt theorem).
\item[{\bf II.}] Consider $(x^ix^j)y^k$ and $x^i(x^jy^k)$.
We have
\[
\begin{array}{lll}
x^i (x^j y^k) &= &x^i (y^k x^j + y^s A^{jk}_s)\\ &= & (y^k x^i + y^r
A^{ik}_r)x^j + (y^s x^i + y^r A^{is}_r) A^{jk}_s\\ &= & y^k x^i x^j +
y^r (A^{ik}_r x^j + A^{jk}_r x^i + A^{is}_r A^{jk}_s)\\ &=& y^k(x^jx^i
+ c^{ij}_sx^s) + y^r (A^{ik}_r x^j + A^{jk}_r x^i + A^{is}_r A^{jk}_s)
\end{array}
\]
and
\[
\begin{array}{lll}
(x^ix^j)y^k&=& (x^jx^i + c^{ij}_sx^s)y^k\\ &=& y^k x^j x^i + y^r
(A^{jk}_r x^i + A^{ik}_r x^j + A^{js}_rA^{ik}_s) +c^{ij}_s (y^k x^s +
y^r A^{sk}_r)
\end{array}
\]
from which it follows that we need
\begin{equation}\label{2.4}
A^{jk}_s A^{is}_r - A^{ik}_s A^{js}_r = c^{ij}_s A^{sk}_r
\end{equation}
\item[{\bf III.}] Consider $(x^iy^j)y^k$ and $x^i(y^jy^k)$. We have
\[
\begin{array}{lll}
(x^iy^j)y^k &=& (y^jx^i + y^s A^{ij}_s)y^k\\
&=& y^j(y^kx^i+y^r A^{ik}_r)- A^{ij}_sy^ky^s\\
&=& -y^ky^jx^i + A^{ik}_ry^jy^r - A^{ij}_sy^ky^s
\end{array}\]
and
\[
\begin{array}{lll}
x^i(y^jy^k)&=& -x^i(y^ky^j)\\
&=& -(y^kx^i + y^r A^{ik}_r)y^j\\
&=& -y^k(y^jx^i + y^s A^{ij}_s) + A^{ik}_ry^jy^r
\end{array}\]
{}From this it is clear that $(x^iy^j)y^k$ reduces to the same normal
form as $x^i(y^jy^k)$.
\item[{\bf IV.}] Consider $(y^iy^j)y^k$ and $y^i(y^jy^k)$. Both reduce
to $-y^ky^jy^i$, so no additional conditions arise.
\end{enumerate}

To summarize we derived that the monomials in definition 2.1 form
indeed a basis of \E\ if (and only if) the equations (2.4) are
satisfied.

Till now we didn't take care of d. Let us describe how this can be
done. Consider the tensor algebra on the letters $x^i$ and $y^j$,
$i,j\in I$. A basis are the monomials
\[ z^{i_1}z^{i_2}\dots z^{i_m} \]
where $i_1,i_2,\ldots , i_m \in I$ (not necessarily ordered) and $z=x$
or $z=y$. On the tensor algebra it is easy to define d, namely
\[
\mbox{d}(z^{i_1}z^{i_2}\dots z^{i_m}) =
\sum_{s=1}^m (-1)^{d_s-1} z^{i_1}z^{i_2}\dots \mbox{d}(z^{i_s})\dots
z^{i_m}
\]
Here $d_s$ is the number of times that $z=y$ in $z^{i_1}z^{i_2}\dots
z^{i_{s-1}}$. Of course d$(x^i)=y^i$ and d$(y^j)=0$. It is not
difficult to prove that d satisfies the requirements in definition 2.1
(b). Now we want to define d on \E. For this it is sufficient to prove
that d preserves the relations of (\ref{2.1}), (\ref{2.2}) and (2.3),
since in that case d leaves the ideal generated by (\ref{2.1}),
(\ref{2.2}) and (2.3) invariant, thanks to d$(\omega \eta)=
\mbox{d}(\omega)\eta + (-1)^n\omega \mbox{d}\eta$ in the tensor
algebra. Now applying d to the relations (\ref{2.2}) yields the
relations (2.3) (like before), and applying d to the relations (2.3)
yields 0. Hence only the relations (\ref{2.1}) remain.  We have
\[
\begin{array}{lll}
\mbox{d}(x^ix^j-x^jx^i)& = &y^ix^j+x^iy^j-y^jx^i-x^jy^i\\
&=& y^ix^j+ (y^jx^i + A^{ij}_ky^k) -y^jx^i- (y^ix^j + A^{ji}_ky^k)\\
&=& (A^{ij}_k - A^{ji}_k)y^k = c^{ij}_ky^k
\end{array} \]
So we find
\begin{equation}\label{2.3}
A^{ij}_k - A^{ji}_k = c^{ij}_k
\end{equation}

In the next section we will study these equations (2.4) and (2.5) more
closely. In particular, we will show that they may be interpreted Lie
algebraically in a very interesting way.

\section{Lie algebraic interpretation of the equations for a
differential calculus}
\setcounter{equation}{0}

The conditions to be satisfied in order to obtain a consistent
differential calculus are the equations (\ref{2.4}) and (\ref{2.3}).

Denote by $\g$ the Lie algebra with basis $x^i (i\in I)$ and with
structure constants $c^{ij}_k$. Now define a linear mapping $\rho : \g
\rightarrow gl(\g)$ by
\[
\rho(x^i) x^a = A^{ia}_b x^b
\]
{}From (\ref{2.4}) it follows that
\[
(A^{ia}_b A^{jb}_c - A^{ja}_b A^{ib}_c)x^c = -c^{ij}_k A^{ka}_c x^c
\]
and thus
\[
A^{ia}_b \rho(x^j) x^b - A^{ja}_b \rho (x^i) x^b = -c^{ij}_k
\rho(x^k)x^a
\]
Therefore
\[
\rho (x^j)A^{ia}_b x^b -\rho (x^i)A^{ja}_b x^b = -\rho (c^{ij}_k
x^k)x^a
\]
yielding
\[
\rho(x^j)\rho (x^i)x^a - \rho(x^i)\rho(x^j)x^a = \rho(c^{ji}_k x^k)x^a
=\rho ([x^j, x^i])x^a
\]
So
\[
\rho (x^j)\rho (x^i) -\rho(x^i) \rho(x^j) = \rho ([x^j, x^i])
\]
and therefore in general we have
\begin{equation}
\rho(x)\rho(y) -\rho(y)\rho(x) = \rho([x,y])
\label{3.1}
\end{equation}
which means that the linear mapping $\rho$ is a representation
of the Lie algebra $\g$ on itself.

In fact it is clear that $\rho$ should be a representation from the
beginning. Namely consider the linear span $Y =
\langle y^i ; i\in I \rangle$. Then $\rho$ can be considered to be a
mapping $\rho : \g\ \to gl(Y)$, given by
\[ \rho(x)y = [x,y] \qquad x\in \g, y\in Y \]
as $x^iy^j-y^jx^i = A^{ij}_k y^k$. Clearly  we have
\[ [[x,\tilde{x}],y] = [x,[\tilde{x},y]] - [\tilde{x},[x,y]] \]
and consequently
\[
\rho([x,\tilde{x}]) = \rho(x)\rho(\tilde{x})-\rho(\tilde{x})\rho(x).
\]

Next we consider the other condition. From the equation (\ref{2.3}) it
follows that
\[
\A x^k - A^{ji}_k x^k = c^{ij}_k x^k
\]
yielding
\[
\rho (x^i)x^j -\rho(x^j)x^i = [x^i, x^j]
\]
and therefore in general
\begin{equation}
\rho(x)y - \rho(y)x = [x,y]
\label{3.2}
\end{equation}

Thus we see that in order to determine a consistent differential
calculus we should look for representations $\rho$ of the given Lie
algebra on itself which moreover satisfy the condition (\ref{3.2}).

This interpretation allows us to describe all constructions till now
in a coordinate-independent way. Let us start with a linear space \g.
Then we found that the differential calculi above are completely
described by two mappings, $c:\g\times \g \to \g$ and $\rho: \g \to
gl(\g)$ such that for all $x,y,z \in \g$:
\begin{enumerate}
\item $c(x,y)=-c(y,x) \quad \mbox{and}\quad
c(c(x,y),z)+c(c(y,z),x)+c(c(z,x),y)=0 $\\
i.e. \g\ is a Lie algebra,
\item $\rho(c(x,y))=\rho(x)\rho(y)-\rho(y)\rho(x)
\quad\mbox{and}\quad
 \rho(x)y - \rho(y)x = c(x,y)$
\end{enumerate}
After choosing coordinates $(x^i)$ we have structure constants
$c^{ij}_k$ instead of $c$ and the numbers $A^{ij}_k$ instead of
$\rho$. In a different basis $(\bar{x}^i)$ the structure constants
will be different, say $\bar{c}^{ij}_k$, and similarly
$\bar{A}^{ij}_k$. If $\tau$ is the change of coordinates, i.e.
$\bar{x}^i=\tau^i_kx^k$ and $\sigma$ the inverse, then we have
\begin{equation}\label{coord}
\bar{c}^{ij}_k = \tau_s^i\tau_r^j\sigma_k^tc^{sr}_t
\quad\mbox{and}\quad
\bar{A}^{ij}_k = \tau_s^i\tau_r^j\sigma_k^tA^{sr}_t
\end{equation}
A special case arises if $\tau$ is an automorphism of the Lie algebra
\g. Then $\bar{c}^{ij}_k=c^{ij}_k$, while not necessarily
$\bar{A}^{ij}_k= A^{ij}_k$. So if $(c^{ij}_k, A^{ij}_k)$ is a
solution, so is $(c^{ij}_k, \bar{A}^{ij}_k)$, with $\bar{A}^{ij}_k$
given by (\ref{coord}) and $\tau$ an automorphism. We will call those
solutions equivalent and will content ourselves determining one
solution from each equivalence class.

The transition from $A^{ij}_k$ to $\bar{A}^{ij}_k$ can be described on
the level of $\rho$ quite easily. One can calculate directly that
\begin{equation}\label{auto}
\bar{\rho}(x) = \tau^{-1}\rho(\tau x) \tau
\end{equation}
Hence we will be interested in all solutions $\rho$ up to an
automorphism in the sense of equation (\ref{auto}).

\section{Semisimple Lie algebras and $gl_n$}
\setcounter{equation}{0}

In this section we will prove that there exists no differential
calculus of the form described in section 2 on the universal
enveloping algebra of any semisimple Lie algebra \g.

To this end we have to prove that there exists no representation
$\rho$ of \g\ on itself, which moreover satisfies
\begin{equation}\label{mult_rep}
[x,y] = \rho(x)y - \rho(y)x \qquad \mbox{ for all } x,y\in \g.
\end{equation}
Let $\rho$ be any representation of \g\ on itself. One can interpret
equation (\ref{mult_rep}) in cohomological terms by saying that the
identity $i : \g \to \g$ is a cocycle in the complex $C^*(\g;\rho)$.
According to the first Whitehead lemma, we know that $H^1(\g;\rho)=0$
for any semisimple \g\ and representation $\rho$. Therefore $i$ is a
coboundary, i.e. there exists an element $a\in C^0(\g;\rho)=\g$ such
that d$(a)=i$. Hence it follows that
\begin{equation}\label{bound}
i(x) = \rho(x)a, \qquad \mbox{or } x=\rho(x)a
\end{equation}
Now set $y=a$ in equation (\ref{mult_rep}) and substitute
$\rho(x)a=x$. We arrive at $[x,a] = x - \rho(a) x$ or equivalently
\begin{equation}\label{rho_a}
\rho(a) = i + \mbox{ad } a
\end{equation}
Here ad denotes the adjoint representation. Taking traces of the
mappings above, we see that tr$(\rho(a)) = \mbox{dim}(\g)$, as
tr$(\mbox{ad }a) =0$. Hence $a$ is represented by a matrix $\rho(a)$
which is not traceless. This is impossible for semisimple Lie
algebras.

For $gl_n$ things turn out differently. Indeed in this case one can
find several solutions, of which we describe one. One can take
\[ \rho(x)y=xy \]
where in the right-hand side $xy$ denotes the multiplication of $x$
and $y$ as $n\times n$-matrices. It is obvious that this is a
representation satisfying $[x,y] = \rho(x)y - \rho(y)x$.  Details for
$n=2$ may be found in \cite{MP}.

\section{Abelian Lie algebras}
\setcounter{equation}{0}

In this section we suppose that the Lie algebra \g\ is Abelian of
dimension $n$.  The equations read
\begin{equation}\label{4.1}
[\rho(x^i), \rho(x^j)]=0 \qquad \mbox{and} \qquad
\rho(x^i)x^j=\rho(x^j)x^i
\end{equation}
Hence the matrices $\rho(x^i)$ form a system of commuting matrices
satisfying an additional condition. In low dimensional cases $(\dim \g
\leq 6)$ maximal commutative matrix algebras have been determined
explicitly (see e.g. \cite{ST}).

Thus for small dimensions all differential calculi of the form
(\ref{2.2}) on universal enveloping algebras of Abelian Lie algebras,
actually affine spaces, can be described explicitly. For dimensions
$n=1,2$ this has already been done by Dimakis et al
\cite{DMS}. See also \cite{BDM}.

In general we can describe the ``regular'' case. By this we mean that
all matrices are simultaneous diagonalizable and linearly independent.
In this case $\rho(x^i) = E^i_i$ with respect to the basis
$\{x^1,\ldots, x^n\}$. Here $E^i_j$ denotes the matrix with the
(j,i)-entry equal to 1, and all others 0. The corresponding calculus
is defined by
\[  x^iy^j = y^jx^i + \delta_i^j y^j. \]

\section{The 2-dimensional solvable Lie algebra}
\setcounter{equation}{0}

Here we are concerned with the smallest non-Abelian (solvable) Lie
algebra \g, which is up to an isomorphism given by the structure
constants $c^{12}_2 = 2$. We will denote $x^1=x$ and $x^2=y$, and we
have $[x,y]=2y$. Our task is to find representations $\rho$ of \g\ on
\g\ itself, satisfying $\rho(x)y-\rho(y)x = 2y$. We will perform these
computations in some detail. For the time being we forget the second
condition; we calculate $\rho(x)$ and $\rho(y)$ with respect to a
suitable basis. These ``standard forms'' can be improved by applying
automorphisms of \g. We are just interested to determine one element
in each orbit $\bar{\rho}(g)= \tau^{-1}\rho(\tau g)\tau$, as we
explained in section 3. We see that instead of $\rho(g)$ we can take
$\rho(\tau g)$ as long as we are only interested in a standard form
for $\rho(g)$ with respect to some basis; indeed the conjugation with
$\tau$ can be seen as merely a basis-transformation.

With respect to some basis, we can assume that $\rho(y)$ has Jordan
normal form. From $[x,y]=2y$ we see that tr$(\rho(y))=0$. Hence we
have two cases, case A and case B.\\[\baselineskip]
{\bf Case A.}
\[ \rho(y)=\left[\begin{array}{cc}
                        \alpha&0\\
                         0&-\alpha
                 \end{array}
           \right]
\]
{}From $[\rho(x),\rho(y)]=2\rho(y)$ it follows that $\alpha=0$, so that
$\rho(y)=0$. So now we can put $\rho(x)$ in Jordan normal form.  We
have two subcases, A1 and A2.\\ {\bf Case A1.}
\[ \rho(y)=0 \qquad \mbox{ and } \qquad
   \rho(x)=\left[\begin{array}{cc}
                        \alpha & 0\\
                         0 & \beta
                 \end{array}
           \right]
\]
Now the second condition comes in. We put
\[
y=\left[\begin{array}{c} y_1\\ y_2 \end{array}
  \right]
\qquad \mbox{ and }\qquad
x=\left[\begin{array}{c} x_1\\ x_2 \end{array}
  \right]
\]
then $\rho(x)y-\rho(y)x=2y$ yields $y_1(\alpha-2)=0 \mbox{ and }
y_2(\beta-2)=0$.  Working this out (splitting with respect to $y_1$
and $y_2$ vanishing or not) yields one solution
\begin{equation}\label{A1}
\rho(y)=0 \qquad \mbox{ and }\qquad
   \rho(x)=\left(\begin{array}{cc}
                        \alpha & 0\\
                         \beta(\alpha-2) & 2
                 \end{array}
           \right)
\end{equation}
Here and hereafter the matrices with round brackets ( and ) are with
respect to the basis $\{x,y\}$.  Finally we can apply the automorphism
$\tau$, given by
\[ \tau=\left(\begin{array}{cc}
                        1 & 0\\
                         \beta & 1
                 \end{array}
           \right)
\]
so that $\bar{\rho}(x)=\tau^{-1}\rho(\tau x)\tau$ takes the same form
as in (\ref{A1}) but with $\beta=0$.
\\[\baselineskip]
{\bf Case A2.}
\[ \rho(y)=0 \qquad \mbox{ and } \qquad
   \rho(x)=\left[\begin{array}{cc}
                        \alpha & 1\\
                         0 & \alpha
                 \end{array}
           \right]
\]
Now $\rho(x)y-\rho(y)x=2y$ yields $y_2=-y_1(\alpha-2) \mbox{ and }
y_2(\alpha-2)=0$.  It follows that $y_1\not=0$ since otherwise $y=0$.
So $y_2(\alpha-2)= -y_1(\alpha-2)^2=0$ yields $\alpha=2$. Hence
$y_2=0$. This leads to a solution of the form
\begin{equation}\label{A2}
\rho(y)=0 \qquad \mbox{ and } \qquad
   \rho(x)=\left(\begin{array}{cc}
                        2 & 0\\
                         \beta & 2
                 \end{array}
           \right) \mbox{ with } \beta={x_2 \over y_1} \not=0
\end{equation}
Applying $\tau=\left(\begin{array}{cc}
                        1 & 0\\
                         0 & {1 \over \beta}
                 \end{array}
           \right)$
yields the same form as in (\ref{A2}) but now with $\beta=1$.
\\[\baselineskip]
{\bf Case B.}
\[ \rho(y)=\left[\begin{array}{cc}
                        0&1\\
                         0&0
                 \end{array}
           \right]
\mbox{ and therefore }
\rho(x)=\left[\begin{array}{cc}
                        \alpha+1 & \beta\\
                         0 & \alpha-1
                 \end{array}
           \right]
\]
Applying $\tau=\left(\begin{array}{cc} 1 & 0\\ -\beta & 1 \end{array}
\right)$ we see that we can take $\beta=0$. The condition
$\rho(x)y-\rho(y)x=2y$ yields is this case $x_2=(\alpha-1)y_1$ and
$y_2(\alpha-3)=0$. Hence case B splits in two cases\\[\baselineskip]
{\bf Case B1.} $y_2=0$ and hence $y_1\not=0$ and $\alpha\not=1$.\\ We
obtain
\begin{equation}\label{B1}
    \rho(y)=\left(\begin{array}{cc}
                        0&0\\
                         \alpha-1&0
                 \end{array}
           \right)
\qquad \mbox{ and  }\qquad
\rho(x)=\left(\begin{array}{cc}
                        \alpha-1 & 0\\
                         \beta & \alpha+1
                 \end{array}
           \right)
\end{equation}
with $\beta=2{x_1 \over y_1}$. For $\alpha\not=-1$ we can apply
$\tau=\left(\begin{array}{cc}
                        1&0\\
                         \gamma&1
                 \end{array}
           \right)$
 with $\gamma=-{\beta \over \alpha+1}$ to get $\rho(x)$ in the
same form as in (\ref{B1}) with $\beta=0$. However for $\alpha=-1$ we
have a bifurcation. Either $\beta=0$ and we have the same form as
above, or $\beta\not=0$. In this case one can scale such that
$\beta=1$. So we have an extra solution here:
\begin{equation}\label{B1a}
    \rho(y)=\left(\begin{array}{cc}
                        0&0\\
                         -2&0
                 \end{array}
           \right)
\qquad \mbox{ and }\qquad
\rho(x)=\left(\begin{array}{cc}
                        -2 & 0\\
                         1 & 0
                 \end{array}
           \right)
\end{equation}
Finally we come to the last case.\\[\baselineskip] {\bf Case B2.}
$y_2\not=0$ and hence $\alpha=3$.\\
This case yields a 2-parameter family of solutions, which by a
automorphism can be put in the form
\begin{equation}\label{B2}
    \rho(y)=\left(\begin{array}{cc}
                        0&1\\
                         0&0
                 \end{array}
           \right)
\qquad \mbox{ and }\qquad
\rho(x)=\left(\begin{array}{cc}
                        4 & 0\\
                         0 & 2
                 \end{array}
           \right)
\end{equation}
Concluding we have 5 families of solutions up to automorphisms; {\it
i}: (\ref{A1}) with $\beta=0$, {\it ii}: (\ref{A2}) with $\beta=1$,
{\it iii}: (\ref{B1}) with $\beta=0$, {\it iv}: (\ref{B1a}) and {\it
v}: (\ref{B2}).

To prove that these classes are indeed not equivalent, one can use
fruitfully the invariance of rank$(\rho(y))$ and tr$(\rho(x))$ under
automorphisms.

\section{The Heisenberg algebra}
\setcounter{equation}{0}

In this section we are concerned with the Heisenberg algebra generated
by $2n+1$ generators $p_i, q_i \ (i=1,\ldots,n)$ and $c$ subjected to
the relations
\begin{equation}
[p_i,q_i]=c \quad (i=1,\ldots,n)
\label{6.1}
\end{equation}
And all other commutators are 0.

We will describe the calculi for $n=1$ in detail. We will write
$x^1=c, x^2=p_1$ and $x^3=q_1$ in this case.  Thus $c_1^{23} =1,
c^{32}_1 = -1$ and all other structure constants vanish. The
calculation of all solutions $\rho$ up to automorphisms can be
performed in the same way as in section 6. The result can be divided
into 4 groups of which it is clear that they are mutually not
equivalent under automorphisms; we give $\rho$ with respect to the
basis $\{c,p,q\}$.
\begin{itemize}
\item[I.] $\rho(c) \not=0$. Then (with $\epsilon=0$ or $\epsilon=1$)
\[ \rho(c)=\left(\begin{array}{ccc} 0 & 0 & 1\\
                                0 & 0 & 0\\
                                0 & 0 & 0
             \end{array}\right) \quad ; \quad
\rho(p)=\left(\begin{array}{ccc} 0 & \epsilon & 0\\
                                0 & 0 & 1\\
                                0 & 0 & 0
             \end{array}\right) \quad\mbox{and}\quad
\rho(q)=\left(\begin{array}{ccc} 1 & -1 & 0\\
                                0 & 1 & 0\\
                                0 & 0 & 1
             \end{array}\right) \]
\item[II.] $\rho(c) =0$ and $\rho$ is not traceless. Then (with
$\epsilon=0$ or $\epsilon=1$)
\[ \rho(p)=\left(\begin{array}{ccc} 0 & \epsilon & 0\\
                                0 & 0 & 0\\
                                0 & 0 & 0
             \end{array}\right) \quad\mbox{and}\quad
\rho(q)=\left(\begin{array}{ccc} 0 & -1 & 0\\
                                0 & 0 & 0\\
                                0 & 0 & 1
             \end{array}\right) \]
\item[III.] $\rho(c) =0$ and $\rho$ is traceless. Then (with $\alpha
\in \C)$ either
\[ \rho(p)=\left(\begin{array}{ccc} 0 & 1 & 0\\
                                0 & 0 & 0\\
                                0 & 0 & 0
             \end{array}\right) \quad\mbox{and}\quad
\rho(q)=\left(\begin{array}{ccc} 0 & -1 & \alpha\\
                                0 & 0 & 0\\
                                0 & 0 & 0
             \end{array}\right) \]
or
\[ \rho(p)=\left(\begin{array}{ccc} 0 & 0 & 1+\alpha\\
                                0 & 0 & 0\\
                                0 & 0 & 0
             \end{array}\right) \quad\mbox{and}\quad
\rho(q)=\left(\begin{array}{ccc} 0 & \alpha & 0\\
                                0 & 0 & 1\\
                                0 & 0 & 0
             \end{array}\right) \]
\item[IV.] $\rho(c) =0$ and the solution is invariant under automorphisms.
\[ \rho(p)=\left(\begin{array}{ccc} 0 & 0 & {1 \over 2}\\
                                0 & 0 & 0\\ 0 & 0 & 0
\end{array}\right) \quad\mbox{and}\quad
\rho(q)=\left(\begin{array}{ccc} 0 & -{1 \over 2} & 0\\
                                0 & 0 & 0\\
                                0 & 0 & 0
             \end{array}\right) \]
\end{itemize}
Remarkably the solution in IV is the only solution invariant under
automorphisms. It corresponds to taking $\A = \frac{1}{2}c^{ij}_k$.
Indeed it is a solution, which can be generalized to higher
dimensional Heisenberg algebras by the same relation to $c^{ij}_k$;
from equation (\ref{coord}) it is clear that also in higher dimensions
it is preserved by automorphisms.

The corresponding consistent differential calculus satisfies all
commutation relations of classical calculus with the exception that
\begin{equation}
\begin{array}{lcl}
p_i\:\mbox{d}q_i &=& \mbox{d}q_i\: p_i + \frac{1}{2} \mbox{d}c\\[3mm]
q_i\: \mbox{d}p_i &=& \mbox{d}p_i\: q_i - \frac{1}{2} \mbox{d}c
\end{array}
\label{6.2}
\end{equation}
Likewise we have for the basic differential operators $\partial_{c},
\partial_{p_i}, \partial_{q_i}$, besides the classical, the non-traditional
commutation relations with the generators
\begin{equation}
\begin{array}{lcl}
\partial_c p_i &=& p_i \partial_c + \frac{1}{2} \partial_{q_i}\\[3mm]
\partial_c q_i &=& q_i \partial_c - \frac{1}{2} \partial_{p_i}
\end{array}
\label{6.3}
\end{equation}

\section{The Witt algebra and the Virasoro algebra}
\setcounter{equation}{0}

We consider the algebra $W$ generated by generators $x^i (i\in \Z)$
subjected to the relations
\begin{equation}
x^i x^j - x^j x^i = (j-i)x^{i+j}
\label{7.1}
\end{equation}
$W$ is called the Witt algebra. We set
\begin{equation}
\A = (j+\mu)\delta^{i+j}_k,
\label{7.2}
\end{equation}
where $\mu$ is a complex number. With this choice $\A$ satisfies
(\ref{2.4}) and (\ref{2.3}). Indeed,
\[
\A - A^{ji}_k = (j+\mu)\delta^{i+j}_k -(i+\mu)\delta^{i+j}_k =
(j-i)\delta_k^{i+j} = c^{ij}_k
\]
and
\[
\begin{array}{lcl}
A^{ik}_r A^{jr}_l - A^{jk}_r A^{ir}_l
&=&(k+\mu)\delta^{i+k}_r (r+\mu)\delta^{j+r}_l - (k+\mu)\delta^{j+k}_r
(r+\mu) \delta^{i+r}_l\\
 &=& (k+\mu)(i+k+\mu)\delta^{i+j+k}_l -
(k+\mu)(j+k+\mu)\delta^{i+j+k}_l\\
&=& (k+\mu)(i-j)\delta^{i+j+k}_l = (i-j)A^{i+j,k}_l = -c^{ij}_s A^{sk}_l
\end{array}
\]
The corresponding calculus satisfies
\begin{equation}
x^i \mbox{d}x^j = \mbox{d}x^j x^i + (j+\mu)\mbox{d}x^{i+j}
\label{7.3}
\end{equation}
and for the basic differential operators $\partial_p (p \in \Z)$ we
have the commutation relations with the generators
\begin{equation}
\partial_p x^k = \delta^k_p + x^k \partial_p + (p-k+\mu)\partial_{p-k}
\label{7.4}
\end{equation}

The Virasoro algebra $V$ is obtained from the algebra $W$ by central
extension with a central element $t$. Then the relations become
\begin{equation}
x^n x^m - x^m x^n = (m-n) x^{n+m} + \frac{1}{12} (m^3 -m)\delta^{m,
-n} t
\label{7.5}
\end{equation}
In case of the Virasoro algebra we put
\begin{equation}
A^{n,m}_k = (m +\mu)\delta^{n+m}_k \quad \mbox{if}\;\;n,m,k \in \Z
\label{7.6}
\end{equation}
and
\[
A^{t,j}_l = A^{i,t}_l = A^{i,j}_t =0
\]
with the exception that
\[
A^{m,-m}_t = \frac{1}{24} (m^3 -m).
\]
This modification of $\A$ satisfies again the relations (\ref{2.4})
and (\ref{2.3}) and we obtain a calculus with the property that for
$n,m \in \Z$ we have
\begin{equation}
x^n \mbox{d}x^m = \mbox{d}x^m x^n +(m+\mu) \mbox{d}x^{n+m}
+\frac{1}{24} (m^3 -m) \delta^{n,-m} \mbox{d}t
\label{7.7}
\end{equation}
and for $k,p \in \Z$
\[
\begin{array}{lcl}
\partial_p x^k &=& \delta^k_p + x^k \partial_p +(p
-k+\mu)\partial_{p-k}\\
\partial_p t &=& t\partial_p\\
\partial_t x^k &=& x^k \partial_t +\frac{1}{24} (k^3
-k)\partial_{-k}\\
\partial_t t &=& 1+t \partial_t
\end{array}
\]
We hope to report on applications to differential equations in a
future publication.

\end{document}